\documentclass[12pt,a4paper]{article}
\usepackage[english]{babel}
\usepackage{amssymb}
\usepackage{inputenc}
\usepackage{amssymb}
\begin{document}

\author{S. Albeverio$^1$, Sh.A. Ayupov$^{2,  *},$   R.Z. Abdullaev$^3$, K.K. Kudaybergenov$^4$}

\title{\bf Additive derivations  on generalized Arens algebras  }

\maketitle
\begin{abstract}
 Given a von Neumann algebra $M$ with  a faithful
normal finite trace $\tau$ denote by $L^\Lambda(M, \tau)$ the
generalized Arens  algebra with respect to $M.$ We give a complete
description of all additive derivations on the algebra
$L^\Lambda(M, \tau).$   In particular each additive derivation on
the algebra $L^{\Lambda}(M, \tau),$ where $M$ is a type II von
Neumann algebra, is inner.
\end{abstract}

\medskip

$^1$ Institut fÄur Angewandte Mathematik, UniversitÄat Bonn,
Endenicherllee. 60, D-53115 Bonn (Germany); SFB 611; HCM; BiBoS;
IZKS; CERFIM (Locarno);  e-mail address:
\emph{albeverio@uni-bonn.de}

$^2$ Institute of Mathematics and Information  Technologies,
 Uzbekistan Academy of Sciences, Dormon Yoli str. 29, 100125,
Tashkent (Uzbekistan), ICTP (Trieste, Italy), e-mail: \emph{sh\_ayupov@mail.ru}

$^3$ Institute of Mathematics and Information Technologies,
Uzbekistan Academy of Science, Dormon Yoli str. 29, 100125,
Tashkent, (Uzbekistan) \emph{arustambay@yandex.ru}

$^4$ Karakalpak state university, Ch. Abdirov str. 1, 142012,
Nukus (Uzbekistan), e-mail: \emph{karim2006@mail.ru}



* Corresponding author

\newpage

\begin{center}
{\bf 1. Introduction}
\end{center}

The present paper continues the series
 of papers \cite{Alb1}-\cite{AK10} devoted to the study and  description of derivations on the
algebra $LS(M)$ of locally measurable  operators affiliated with a
von Neumann algebra $M$ and on its various subalgebras.

Let  $\mathcal{A}$ be an algebra over the field complex number
$\mathbb{C}.$ A linear (additive) operator
$D:\mathcal{A}\rightarrow \mathcal{A}$ is called a linear
(additive) \emph{derivation} if it satisfies the identity
$D(xy)=D(x)y+xD(y)$ for all  $x, y\in \mathcal{A}$ (Leibniz rule).
Each element  $a\in \mathcal{A}$ defines a linear derivation $D_a$
on $\mathcal{A}$ given as $D_a(x)=ax-xa,\,x\in \mathcal{A}.$ Such
derivations $D_a$ are said to be \emph{inner derivations}. If the
element  $a$ implementing the derivation $D_a$ on $\mathcal{A},$
belongs to a larger algebra $\mathcal{B},$ containing
$\mathcal{A}$ (as a proper ideal, as usual) then $D_a$ is called a
\emph{spatial derivation}.

One of the main problems in the theory of derivations is to prove
the automatic continuity, innerness or spatialness of
derivations  or to show the existence of non inner
and discontinuous derivations on various topological algebras.

In this direction  A.~F.~Ber, F.~A.~Sukochev,
V.~I.~Chilin~\cite{Ber} obtained necessary and sufficient
conditions for the existence of non trivial derivations on
commutative regular algebras. In particular they have proved that
the algebra  $L^{0}(0, 1)$ of all (classes of equivalence of)
complex measurable functions on  the interval $(0, 1)$ admits non
trivial derivations. Independently A.~G.~Kusraev~\cite{Kus1} by
means of Boolean-valued analysis has also proved the existence of
non trivial derivations and automorphisms on $L^{0}(0, 1).$ It is
clear that these derivations are discontinuous in the measure
topology, and therefore they are neither inner nor spatial. It was
conjectured that the existence of such exotic examples of
derivations deeply depends on the commutativity of the underlying
von Neumann algebra $M.$ In this connection we have initiated the
study of the above problems in the non commutative case
\cite{Alb1}-\cite{Alb3},
 by considering derivations on the algebra $LS(M)$ of
all locally measurable operators affiliated with a von Neumann
algebra $M$ and on various subalgebras of $LS(M).$ In \cite{Alb1}
 noncommutative Arens algebras $L^\omega(M,
\tau) = \bigcap\limits_{p\geq 1} L^p(M, \tau)$ and related
algebras associated with a von Neumann algebra $M$ and a faithful
normal semi-finite trace $\tau$ have been considered. It has been proved that every derivation on
this algebra is spatial, and, if the trace $\tau$ is finite, then
all derivations are inner. In  \cite{Alb2} and  \cite{Alb3}    the mentioned conjecture concerning
derivations on on the algebra
 $LS(M)$ has been confirmed for type I von Neumann algebras.

  Recently this conjecture was also independently confirmed for the type I
case in the paper of A.F. Ber, B. de Pagter and A.F. Sukochev
\cite{Ber1} by means of a representation of measurable operators
as operator valued functions. Another approach to similar problems
in the framework of type I $AW^{*}$-algebras has been outlined in
the paper of A.F. Gutman, A.G.Kusraev and  S.S. Kutateladze
\cite{Gut}.

In \cite{Alb2}  we considered derivations on the algebra $LS(M)$
of all locally measurable operators affiliated with a type I von
Neumann algebra $M$, and also on its subalgebras $S(M)$ -- of
measurable operators and $S(M, \tau)$ of $\tau$-measurable
operators, where $\tau$ is a faithful normal semi-finite trace on
$M.$ It was proved that an arbitrary derivation $D$ on each of these
algebras can be uniquely decomposed into the sum $D=D_a+D_\delta$
where the derivation $D_a$ is inner (for $LS(M),$ $S(M)$ and $S(M,
\tau)$) while the derivation $D_\delta$ is an extension of a derivation $\delta$ (possibly non trivial) on the center of the corresponding algebra.

In the present paper we consider additive derivations on
generalized Arens  algebras in the sense of Kunze  \cite{Kun} with respect to a von Neumann algebra
with a faithful normal finite trace.

In section 1 we give some necessary properties of the generalized
Arens algebra $L^\Lambda(M, \tau).$

Section 2 is devoted to study of additive derivations on generalized
Arens algebras.   We prove that an arbitrary additive
derivation $D$ on the algebra  $L^\Lambda(M, \tau)$ can be uniquely
decomposed into the sum $D=D_a+D_\delta,$ where the derivation
$D_a$ is inner while the derivation $D_\delta$ is an extension of some
additive derivation $\delta$ on the center of the  algebra
$L^\Lambda(M, \tau).$ In particular, if $M$ is a type II von
Neumann algebra then every additive derivation on the algebra
$L^\Lambda(M, \tau)$ is inner.

\begin{center}
\textbf{2. Generalized Arens  algebras}
\end{center}

Let  $H$ be a complex Hilbert space and let  $B(H)$ be the algebra
of all bounded linear operators on   $H.$ Consider a von Neumann
algebra $M$  in $B(H)$ with the operator norm $\|\cdot\|_M.$ Denote by
$P(M)$ the lattice of projections in $M.$

A linear subspace  $\mathcal{D}$ in  $H$ is said to be
\emph{affiliated} with  $M$ (denoted as  $\mathcal{D}\eta M$), if
$u(\mathcal{D})\subset \mathcal{D}$ for every unitary  $u$ from
the commutant
$$M'=\{y\in B(H):xy=yx, \,\forall x\in M\}$$ of the von Neumann algebra $M.$

A linear operator  $x$ on  $H$ with the domain  $\mathcal{D}(x)$
is said to be \emph{affiliated} with  $M$ (denoted as  $x\eta M$) if
$\mathcal{D}(x)\eta M$ and $u(x(\xi))=x(u(\xi))$
 for all  $\xi\in
\mathcal{D}(x).$

 Let   $\tau$ be a faithful normal semi-finite trace on $M.$ We recall that a closed linear operator
  $x$ is said to be  $\tau$\emph{-measurable} with respect to the von Neumann algebra
   $M,$ if  $x\eta M$ and   $\mathcal{D}(x)$ is
  $\tau$-dense in  $H,$ i.e. $\mathcal{D}(x)\eta M$ and given   $\varepsilon>0$
  there exists a projection   $p\in M$ such that   $p(H)\subset\mathcal{D}(x)$
  and $\tau(p^{\perp})<\varepsilon.$
   The set $S(M,\tau)$ of all   $\tau$-measurable operators with respect to  $M$
   is a unital *-algebra
when equipped with the algebraic operations of strong addition and
multiplication and taking the adjoint of an operator (see
\cite{Mur}).

    Consider the topology  $t_{\tau}$ of convergence in measure or \emph{measure topology}
    on $S(M, \tau),$ which is defined by
 the following neighborhoods of zero:
$$V(\varepsilon, \delta)=\{x\in S(M, \tau): \exists e\in P(M), \tau(e^{\perp})\leq\delta, xe\in
M,  \|xe\|_{M}\leq\varepsilon\},$$  where $\varepsilon, \delta$
are positive numbers, and $\|.\|_{M}$ denotes the operator norm on
$M$.

 It is well-known
\cite{Mur} that $S(M, \tau)$ equipped with the measure topology is
a complete metrizable topological *-algebra.

Recall \cite{Kra}  that $\phi$ is a Young  function, if
\[
\phi(t)=\int\limits_{0}^t\varphi(s)\,ds,\quad t\geq 0,
\]
where the real-valued function $\varphi$ defined on $[0, \infty)$
has the following properties:

(i)  $\varphi(0)=0, \, \varphi(s)>0$ for $s>0$ and
$\lim\limits_{s\rightarrow\infty}\varphi(s)=\infty,$

(ii) $\varphi$ is right continuous,

(iii) $\varphi$ is nondecreasing on $(0, \infty).$

Every Young function is a continuous, convex and strictly
increasing function. For every Young function $\phi$ there is a
complementary Young function $\psi$ given by the density
\[
\psi(t)=\sup\{s: \phi(s)\leq t\}.
\]
The complement of $\psi$ is $\phi$ again. Further a Young function
$\phi$ is said to satisfy the $\Delta_2$-condition, shortly
$\phi\in \Delta_2,$ if there exists a $k>0$ and $T\geq 0$ such
that:
\[
\phi(2t)\leq k\phi(t)
\]
for all $t\geq T.$

Put
\[
K_\phi=\{x\in S(M, \tau): \tau(\phi(|x|))\leq1\}
\]
and
 \[
L_\phi(M, \tau)=\bigcup\limits_{n=1}^\infty nK_\phi.
\]

It is known \cite{Mur1} (see also  \cite{Kun}) that $L_{\phi}(M,
\tau)$ is a Banach space with respect to the norm
$$\|x\|_\phi=\inf\left\{\lambda>0: \frac{1}{\lambda}x\in K_\phi\right\},\quad x\in L_{\phi}(M, \tau).$$

We recall from\cite{Kun}  that $\phi_1\prec\phi_2,$ if there exist
two nonnegative constants $c$ and $T$ such that
$\phi_1(t)\leq\phi_2(c t)$ for all $t\geq T.$
 Let $\Lambda$ be a generating family of Young functions,
i.e. for $\phi_1, \phi_2\in\Lambda$ there is a $\psi\in\Lambda$
with $\phi_1, \phi_2\prec\psi.$
 A generating family $\Lambda$ of Young functions is said to be
 quadratic, if for any $\phi\in\Lambda$ there is a
 $\psi\in\Lambda$ such that the composition of $\phi$ and the
 squaring function as a Young function is smaller than $\psi$ regarding the
 partial order $\prec,$ i.e. there are  $c>0$ and  $T\geq0$
 with $\phi(t^2)\leq\psi(ct)$ for all $t\geq T.$ For a quadratic
 family $\Lambda$ of Young functions we define
 \[
 L^\Lambda(M, \tau)=\bigcap\limits_{\phi\in\Lambda}L_{\phi}(M, \tau).
 \]
On the space  $L^\Lambda(M, \tau)$ one can consider the topology
$t_\Lambda$ generated by the system  of norms $\{\|\cdot\|_{\phi}:
\phi\in \Lambda\}.$

It is known \cite[Proposition 4.1]{Kun} that if $\Lambda$ is a
quadratic family of Young functions, then ($L^\Lambda(M, \tau),
t_\Lambda$) is a complete locally convex *-algebra with jointly
continuous multiplications.

Note that if $\Lambda=\{t^p: p\geq 1\}$ we have  that
\[
L^\Lambda(M, \tau)=L^\omega(M, \tau)=\bigcap\limits_{p\geq 1}L^p(M, \tau).
\]
Non-commutative Arens algebras $L^\omega(M, \tau)$ were introduced
by Inoue \cite{Ino} and their properties were investigated in
\cite{Abd}. Generalized Arens algebras  were introduced
by Kunze  \cite{Kun}.

 Let $\varphi\in \Lambda$ be a Young function. Then there exists
a Young function $\phi\in \Lambda$ and $k>0$ such that
\begin{equation}
 ||x
y||_{\varphi}\leq k||x||_{\phi}||y||_{\phi}
\end{equation}
for all $x, y\in L^{\Lambda}(M, \tau)$ (see \cite{Kun}).

Let us remark that, if $\tau$ is a finite trace, then $t\prec
\phi(t)$ for every Young function, and for any quadratic family
$\Lambda$ of Young functions we obtain that
\begin{equation}
L^{\Lambda}(M, \tau)\subset L^{\omega}(M, \tau).
\end{equation}
Further, if every $\phi\in \Lambda$ satisfies the
$\Delta_2$-condition then
\[
L^{\omega}(M, \tau)\subset L^{\Lambda}(M, \tau).
\]

It is known \cite{Kun} that if $N$ is a von Neumann subalgebra of
$M$ then
\[
L_\phi(N, \tau_N)=S(N, \tau_N)\cap L_\phi(M, \tau),
\]
where $\tau_N$ is the restriction  of the trace $\tau$ onto  $N.$

It should be noted that if $M$ is a finite von Neumann algebra
with a faithful normal semi-finite trace $\tau,$ then the
restriction $\tau_Z$ of the trace $\tau$ onto the center $Z(M)$ of
$M$ is also semi-finite.

Further we shall need the description of the center of the algebra
$L^{\Lambda}(M, \tau)$  for  von Neumann algebras with a faithful
normal finite trace .

\textbf{Proposition 2.1.} \emph{Let   $M$ be a  von Neumann
algebra  with
 a faithful normal finite trace
$\tau$ and with the center $Z(M).$ Then}
\[
Z(L^{\Lambda}(M, \tau))=L^{\Lambda}(Z(M), \tau_Z).
\]

Proof. Using the equality
\[
L_\phi(N, \tau_N)=S(N, \tau_N)\cap L_\phi(M, \tau).
\]
we obtain that
\[
L^{\Lambda}(N, \tau_N)=S(N, \tau_N)\cap L^{\Lambda}(M, \tau).
\]
Hence
\[
L^{\Lambda}(Z(M), \tau_Z)=S(Z(M), \tau_Z)\cap L^{\Lambda}(M,
\tau)=
\]
\[
=Z(S(M, \tau))\cap L^{\Lambda}(M, \tau)=Z(L^{\Lambda}(M, \tau)),
\]
i.e.
\[
Z(L^{\Lambda}(M, \tau))=L^{\Lambda}(Z(M), \tau_Z).
\]
The proof is complete. $\blacksquare$

\begin{center}
 \textbf{3. Derivations on the generalized Arens algebras}
 \end{center}

In this section we give a complete description of all additive
derivations on the algebra $L^\Lambda(M, \tau).$

 Let  $\mathcal{A}$ be an algebra with the center  $Z(\mathcal{A})$ and let
$D:\mathcal{A}\rightarrow \mathcal{A}$ be an additive  derivation.
Given any $x\in \mathcal{A}$ and a central element $a\in
Z(\mathcal{A})$ we have
$$D(a x)=D(a)x+a D(x)$$
and
$$D(x a)=D(x)a +xD(a).$$
Since  $a x=x a$ and $a D(x)=D(x) a,$ it follows that
$D(a)x=xD(a)$ for any $a\in \mathcal{A}.$ This means that $D(a)\in
Z(\mathcal{A}),$ i.e. $D(Z(\mathcal{A}))\subseteq Z(\mathcal{A}).$
Therefore  given any additive derivation $D$ on the algebra
$\mathcal{A}$ we can consider its restriction
$\delta:Z(\mathcal{A})\rightarrow Z(\mathcal{A}).$

We shall need some facts about additive derivations
$\delta:\mathbb{C}\rightarrow\mathbb{C}.$ Every such derivation
vanishes at every algebraic number. On the other hand, if
$\lambda\in\mathbb{C}$ is transcendental then there is a additive
derivation $\delta:\mathbb{C}\rightarrow\mathbb{C}$ which does not
vanish at $\lambda$ (see \cite{Sam}).

Let $M_n(\mathbb{C})$ be the algebra of $n\times n$ matrices over
$\mathbb{C}.$ If $e_{i,j},\,i,j=\overline{1, n},$ are the matrix
units in $M_n(\mathbb{C}),$ then each element $x\in
M_n(\mathbb{C})$ has the form
 $$x=\sum\limits_{i,j=1}^{n}\lambda_{i j}e_{i j},\,\lambda_{i,j}\in \mathbb{C},\,i,j=\overline{1, n}.$$
Let  $\delta:\mathbb{C}\rightarrow \mathbb{C}$ be an additive
derivation. Setting
\begin{equation}
  D_{\delta}\left(\sum\limits_{i,j=1}^{n}\lambda_{i j}e_{i j}\right)=
 \sum\limits_{i,j=1}^{n}\delta(\lambda_{i j})e_{i j}
\end{equation}
 we obtain a well-defined additive  operator
 $D_\delta$ on the algebra $M_n(\mathbb{C}).$ Moreover
 $D_\delta$ is an additive  derivation on the algebra  $M_n(\mathbb{C})$
 and its restriction onto the center of the algebra  $M_n(\mathbb{C})$ coincides with the given $\delta.$

 It is known \cite[Theorem 2.2]{Sem} that if   $M$ be a von Neumann factor  of type
  I$_{n}, n \in \mathbb{N}$
then every additive  derivation  $D$ on the algebra $M$ can be
uniquely represented as a sum
  $$D=D_{a}+D_{\delta ,}$$ where  $D_{a}$ is an inner derivation implemented by an element  $a\in M$
while $D_{\delta} $ is the additive derivation of the form (3)
generated by an additive  derivation $\delta$ on the center of $M$
identified with $\mathbb{C}.$

Note that if $M$   is a  finite-dimensional von Neumann algebra
then $L^{\Lambda}(M, \tau)=M$ for any faithful normal finite trace
$\tau.$

Now let  $M$ be an arbitrary finite-dimensional von Neumann
algebra with the center $Z(M).$ There exist a family of mutually
orthogonal  central projections $\{z_{1}, z_{2},...,z_{k}\}$ from
$M$ with $\bigvee\limits_{i=1}^k z_i=\textbf{1}$ and $n_1,
n_2,..., n_k\in\mathbb{N}$ such that the algebra $M$ is
*-isomorphic with the $C^{*}$-product of von Neumann factors $z_i
M$ of type I$_{n_i}$ respectively, i.e.
\[
M\cong M_{n_1}(\mathbb{C})\oplus
M_{n_2}(\mathbb{C})\oplus...\oplus M_{n_k}(\mathbb{C}). \]
 Suppose
that    $D$ is an additive  derivation on  $M,$ and $\delta$ is
its restriction onto its center  $Z(M).$ Since
$\delta(zx)=z\delta(x)$ for all central projection $z\in Z(M)$ and
$x\in M$ then  $\delta$ maps each $z_i Z(M)\cong \mathbb{C}$ into
itself, $\delta$ generates an additive derivation $\delta_i$ on
$\mathbb{C}$ for each $i=\overline{1, k}.$

Let     $D_{\delta_i}$ be the additive derivation on the matrix
algebra $M_{n_i}(\mathbb{C}), i=\overline{1, k},$ defined as in
(3). Put
\begin{equation}
D_\delta((x_i)_{i=1}^{k})=(D_{\delta_i}(x_i)),\,(x_i)_{i=1}^{k}\in
M.
\end{equation}
 Then the map  $D_\delta$  is an additive  derivation on $M.$

\textbf{Lemma 3.1.} \emph{Let  $M$ be a finite-dimensional  von
Neumann algebra. Each additive derivation  $D$ on the algebra $M$
can be uniquely represented in the form
$$D=D_{a}+D_{\delta ,}$$
where  $D_{a}$ is an inner derivation implemented by an element
$a\in M,$ and $D_{\delta} $ is an additive  derivation given (4).}

Proof. Let  $D$ be an additive  derivation on $M,$ and let
$\delta$ be its restriction onto $Z(M).$ Consider an additive
derivation $D_\delta$ on $Z(M)$ of the form (4), generated by an
additive  derivation $\delta.$ Since  additive  derivations $D$
and $D_\delta$ coincide on $Z(M)$ ,then  an additive  derivation
of the form  $D-D_\delta$ is a linear derivation. Hence  by
Sakai's theorem \cite[Theorem  4.1.6]{Sak1} $D-D_\delta$ is an
inner derivation. This means that there exists an element $a\in M$
such that $D_a=D-D_\delta$ and therefore  $D=D_{a}+D_{\delta}.$
The proof is complete. $\blacksquare$

Now let $M$ be a commutative von Neumann algebra with a faithful
normal finite trace $\tau.$ Given an arbitrary additive derivation
$\delta$ on $L^\Lambda(M, \tau)$ the element
$$z_\delta=\inf\{z\in P(M): z\delta=\delta\}$$
is called the support of the derivation  $\delta.$

Suppose that   $M$ is  a commutative  von Neumann algebra with
 a faithful normal finite trace
$\tau$ and  $q_1, q_2,...,q_k$ are atoms in $M.$ Then
\[
L^{\Lambda}(M, \tau)\cong q_1\mathbb{C}\oplus
q_2\mathbb{C}\oplus...\oplus q_k\mathbb{C}\oplus p L^{\Lambda}(M,
\tau),\] where $p=\textbf{1}-\bigvee\limits_{i=1}^{k}q_i.$

Now if $\delta_i:\mathbb{C}\rightarrow\mathbb{C}$ is an additive
derivation then
\begin{equation}
\delta(x)=(\delta_1(q_1x),...,\delta_k(q_k x), 0),\,\,\,x\in
L^{\Lambda}(M, \tau)
\end{equation} is also  an additive derivation.
Note
 that  $z_\delta=\bigvee\{q_i:   \delta_i\neq 0, 1\leq i\leq k\}.$

\textbf{Lemma  3.2.} \emph{Let   $M$ be a commutative  von Neumann
algebra  with
 a faithful normal finite trace
$\tau.$ For any non trivial additive derivation $\delta:
L^{\Lambda}(M, \tau)\rightarrow L^{\Lambda}(M, \tau)$
 there exists a sequence
$\{a_n\}_{n=1}^{\infty}$ in $M$ with $|a_n|\leq \textbf{1},\,n\in
\mathbb{N}$, such that
$$|\delta(a_n)|\geq n z_\delta$$ for all $n\in \mathbb{N}.$}

In \cite[Lemma  2.6]{Alb2}  (see also \cite[Lemma 4.6]{Ber1}) this
assertion was proved for linear derivations on the algebra $S(M),$
but same the proof is applies also to the case of additive
derivations on $L^{\Lambda}(M, \tau).$

The following result shows that the above construction   (5) is
the general form of additive derivations on the generalized Arens
algebras in the commutative case.

 \textbf{Lemma   3.3.} \emph{Let   $M$ be a
commutative von Neumann algebra  with
 a faithful normal finite trace
$\tau$ and  let $\delta$ be an additive derivation on the algebra
$L^{\Lambda}(M, \tau).$ Then $z_\delta M$ is a finite-dimensional
algebra.}

Proof. Suppose that $z_\delta M$ is infinite-dimensional. Then
there exists an infiniti sequence of mutually orthogonal
projections $\{z_n\}_{n=1}^{\infty}$ in $M$ such that
  $\bigvee\limits_{n=1}^{\infty}z_n=z_\delta.$  By Lemma 3.2 there exists a sequence
$\{a_n\}_{n=1}^{\infty}$ in $M$ with $|a_n|\leq \textbf{1},\,n\in
\mathbb{N}$, such that
\begin{equation}
|\delta(a_n)|\geq 2^n \tau(z_n)^{-1} z_\delta
\end{equation}
for
all $n\in \mathbb{N}.$ Put
\[
a=\sum\limits_{n=1}^\infty\frac{a_n z_n}{2^n}.
\]
Then $a\in M\subset L^{\Lambda}(M, \tau)$ and
\[
\delta(a)=\delta\left(\sum\limits_{n=1}^\infty\frac{a_n
z_n}{2^n}\right)=\sum\limits_{n=1}^\infty\frac{z_n}{2^n}\delta(a_n).
\]
From (6) we obtain that
\[
|\delta(a)|=\sum\limits_{n=1}^\infty\frac{z_n}{2^n}|\delta(a_n)|\geq
\sum\limits_{n=1}^\infty\frac{z_n}{2^n}2^n \tau(z_n)^{-1}
z_\delta,
\]
i.e.
\[
|\delta(a)|\geq \sum\limits_{n=1}^\infty \tau(z_n)^{-1} z_n.
\]
Thus
\[
\tau(|\delta(a)|)\geq \sum\limits_{n=1}^\infty \tau(z_n)^{-1}
\tau(z_n)=\sum\limits_{n=1}^\infty 1=\infty.
\]
This means that $\delta(a)\notin L^{1}(M, \tau).$ Then by (2) we
have that  $\delta(a)\notin L^{\Lambda}(M, \tau).$ This
contradiction implies that $z_\delta M$ is a finite-dimensional
algebra. The proof is complete. $\blacksquare$

Lemma 3.3 implies the following

 \textbf{Corollary   3.1.} \emph{Let   $M$ be a
commutative von Neumann algebra  with
 a faithful normal finite trace
$\tau$ such that the Boolean algebra $P(M)$ of all projections of
$M$ is continuous. Then every additive derivation on the algebra
$L^{\Lambda}(M, \tau)$ is zero.}

Note that the properties of  additive derivations on the algebras
$S(M, \tau)$ and $L^{\Lambda}(M, \tau),$ where $M$ be a
commutative von Neumann algebra  with
 a faithful normal finite trace
$\tau,$ are quite opposite. Indeed, if the Boolean
algebra $P(M)$ is  continuous then the algebra $S(M, \tau)$
admits a non-zero linear, in particular additive, derivation,
(see \cite[Theorem 3.3]{Ber}), whereas the algebra $L^{\Lambda}(M,
\tau)$ in this case does not admit a non-zero additive derivation (see Corollary
3.1).

Now we consider the noncommutative case.

We shall need following result (\cite[Theorem 4.1]{AAK}, see also \cite[Theorem 6.8]{AK10}).

\textbf{Theorem 3.1.} \emph{Let $M$ be a von Neumann algebra with
 a faithful normal finite trace $\tau$.
  If $A\subseteq L^\omega(M, \tau)$ is a solid
*-subalgebra such that $M\subseteq A$,
 then every linear derivation on $A$ is inner.}

 The following theorem is one
of the main results of this paper.

 \textbf{Theorem 3.2.} \emph{Let   $M$ be a type II von Neumann
algebra  with
 a faithful normal finite trace
$\tau.$ Then every additive derivation on the algebra
$L^{\Lambda}(M, \tau)$ is inner.}

The proof of the theorem 3.2 follows from Theorem 3.1 and the
following assertion.

\textbf{Lemma  3.4.} \emph{Let   $M$ be a type II von Neumann
algebra  with
 a faithful normal finite trace
$\tau,$ and suppose that $D:L^{\Lambda}(M, \tau)\rightarrow
L^{\Lambda}(M, \tau)$ is an additive derivation. Then
$D|_{Z(L^{\Lambda}(M, \tau))}\equiv 0,$ in particular, $D$ is a
linear.}

Proof. Let  $D$ be an additive  derivation on $L^{\Lambda}(M,
\tau),$ and let $\delta$ be its restriction onto $Z(L^{\Lambda}(M,
\tau)).$

Since  $M$ is of type II there exists a sequence of mutually
orthogonal projections $\{p_n\}_{n=1}^{\infty}$ in $M$  with
central covers $\textbf{1}$ (i.e.the $\{p_n\}$ are faithful
projections). For any bounded sequence
$B=\{b_n\}_{n\in\mathbb{N}}$ in $Z(M)$ define an operator $x_B$ by
$$x_B=
\sum\limits_{n=1}^{\infty}b_n p_n.$$ Then
\begin{equation}
x_B p_n=p_n x_B=b_n p_n
\end{equation} for all $n\in\mathbb{N}.$

Take  $b\in Z(M)$ and $n\in \mathbb{N}.$ From the identity
$$D(b p_n)=D(b)p_n+b D(p_n)$$
 multiplying it by $p_n$ on both sides we obtain
$$p_nD(b p_n)p_n=p_n D(b)p_n+b p_n D(p_n)p_n.$$
Since  $p_n$ is a projection, one has that  $p_n D(p_n)p_n=0,$ and
since $D(b)=\delta(b)\in Z(L^{\Lambda}(M, \tau)),$ we have
\begin{equation}
p_nD(b p_n)p_n=\delta(b)p_n.
\end{equation}

Now from the identity
$$D(x_B p_n)=D(x_B)p_n+x_B D(p_n),$$
in view of (7)  one has similarly
$$p_nD(b_n p_n)p_n=p_n
D(x_B)p_n+b_n p_n D(p_n)p_n,$$ i.e.
\begin{equation}
p_n D(b_n p_n)p_n=p_nD(x_B)p_n.
\end{equation} Now (8) and
(9) imply
 \begin{equation}
p_nD(x_B)p_n=\delta(b_n)p_n.
\end{equation}

Let $\varphi\in \Lambda.$ By (1) there are $\phi, \psi\in \Lambda$
and $k>0$ such that
\[
||x_1x_2x_3||_{\varphi}\leq
k||x_1||_{\phi}||x_2||_{\phi}||x_3||_{\psi}
\]
for all $x_1, x_2, x_3\in L^{\Lambda}(M, \tau).$ If we suppose
that $\delta\neq 0$ then  $z_\delta\neq0.$ By Lemma 3.2 there
exists a bounded sequence $B=\{b_n\}_{n\in\mathbb{N}}$ in $Z(M)$
such that
\[
|\delta(b_n)|\geq n c_n z_\delta \]
 for all $n\in \mathbb{N},$ where $c_n=k ||p_n||_{\phi}^{2}||p_n z_\delta||_{\varphi}^{-1}.$
Then in view of (10) we obtain
\[
k||p_n||_{\phi}||D(x)||_{\psi}||p_n||_{\phi}\geq||p_n
D(x)p_n||_{\varphi}= \]
\[ =||\delta(b_n)p_n||_{\varphi}\geq ||n
c_n p_n z_\delta||_{\varphi}=n c_n||p_n z_\delta||_{\varphi},
\]
i.e.
\[
||D(x)||_{\psi}\geq n c_n k^{-1}||p_n||_{\phi}^{-2}||p_n
z_\delta||_{\varphi}.
\]
Hence
\[
||D(x)||_{\psi}\geq n
\]
for all $n\in\mathbb{N}.$
 This contradiction implies that
  $\delta\equiv 0,$ i.e. $D$ is identically
zero on the center of $L^{\Lambda}(M, \tau),$ and therefore it is
linear. The proof is complete. $\blacksquare$

Now consider an additive  derivation  $D$ on  $L^{\Lambda}(M,
\tau)$ and let $\delta$ be its restriction onto its center
$Z(L^{\Lambda}(M, \tau)).$ By Lemma 3.3
 $z_\delta M$ is a finite-dimensional  and  $z_\delta^{\perp}\delta\equiv 0,$ i.e. $\delta=z_\delta\delta.$

Let   $D_\delta$ be the derivation on  $z_\delta L^{\Lambda}(M,
\tau)=z_\delta M$ defined as in (4) and consider its extension
$D_\delta$ on $L^{\Lambda}(M, \tau)=z_\delta L^{\Lambda}(M,
\tau)\oplus z_\delta^{\perp}L^{\Lambda}(M, \tau)$ which is defined
as
\begin{equation}
D_\delta(x_1+x_2):=D_\delta(x_1),\,x_1\in z_\delta L^{\Lambda}(M,
\tau),x_2\in z_\delta^{\perp}L^{\Lambda}(M, \tau).
\end{equation}

The following theorem is the main result of this paper, and gives
the general form of derivations on the algebra $L^{\Lambda}(M,
\tau).$

\textbf{Theorem  3.3.} \emph{Let  $M$ be a  von Neumann algebra
with a faithful normal finite trace $\tau.$ Each additive
derivation $D$ on $L^{\Lambda}(M, \tau)$ can be uniquely
represented in the form
$$D=D_{a}+D_{\delta}$$
where  $D_{a}$ is an inner derivation implemented by an element
$a\in L^{\Lambda}(M, \tau),$ and $D_{\delta} $ is an additive
derivation of the form (11), generated by an additive  derivation
$\delta$ on the center of $L^{\Lambda}(M, \tau)$.}

Proof. Let  $D$ be an additive  derivation on $L^{\Lambda}(M,
\tau),$ and let $\delta$ be its restriction onto $Z(L^{\Lambda}(M,
\tau))=L^{\Lambda}(Z(M), \tau_Z)).$ By Lemma 3.3
 $z_\delta Z(M)$ is  finite-dimensional.
Thus $z_\delta M$ is a  $C^*$-product of a finite number of von
Neumann factors of type I$_n$ or II. Since by Lemma 3.4 any
additive derivation on $L^{\Lambda}(M, \tau),$ where $M$ is a type
II algebra, is linear, then
 by Theorem
3.2 it is inner. Therefore $z_\delta M$ is a $C^*$-product of a
finite number of von Neumann factors of type I$_n$.

Now consider an additive derivation $D_\delta$ on $L^{\Lambda}(M,
\tau)$ of the form (11), generated by a derivation $\delta.$ Since
the derivations $D$ and $D_\delta$ coincide on $L^{\Lambda}(Z(M),
\tau))$ then  $D-D_\delta$ is a linear derivation. Hence Theorem
3.2 implies that the derivation $D-D_\delta$ is  inner. This means
that there exists an element $a\in L^{\Lambda}(M, \tau)$ such that
$D_a=D-D_\delta$ and therefore $D=D_{a}+D_{\delta}.$ The proof is
complete. $\blacksquare$

Theorem 3.3 implies that following.

\textbf{Corollary 3.2.} \emph{Let  $M$ be a  von Neumann algebra
without type I$_n$  direct summands and with a faithful
normal finite trace $\tau.$  Then each additive  derivation  on
$L^{\Lambda}(M, \tau)$ is inner.}

\newpage

\medskip

\textbf{Acknowledgments.} \emph{The third named author would like
to acknowledge the hospitality of the "Institut f\"{u}r Angewandte
Mathematik", Universit\"{a}t Bonn (Germany). This work is supported
in part by the German Academic Exchange Service -- DAAD .}

\end{document}